\documentclass[12pt]{article}
\usepackage{amsmath,amssymb,amsbsy,amsfonts,amsthm,
  latexsym,amsopn,amstext,amsxtra,euscript,amscd}

\begin{document}

\newtheorem{prop}{Proposition}[section]
\newtheorem{lem}[prop]{Lemma}
\newtheorem{cor}[prop]{Corollary}
\newtheorem{defi}[prop]{Definition}
\newtheorem{thm}[prop]{Theorem}
\newtheorem{rem}[prop]{Remark}
\newtheorem{rems}[prop]{Remarks}
\newtheorem{fac}[prop]{Fact}
\newtheorem{facs}[prop]{Facts}
\newtheorem{com}[prop]{Comments}
\newtheorem{prob}{Problem}
\newtheorem{problem}[prob]{Problem}
\newtheorem{ques}{Question}
\newtheorem{question}[ques]{Question}


\def\scr{\scriptstyle}
\def\\{\cr}
\def\({\left(}
\def\){\right)}
\def\[{\left[}
\def\]{\right]}
\def\<{\langle}
\def\>{\rangle}
\def\fl#1{\left\lfloor#1\right\rfloor}
\def\rf#1{\left\lceil#1\right\rceil}
\def\defn{\noindent{\bf Definition\/}\ \ }
\def\rem{\noindent{\bf Remarks\/}\ \ }
\def\card#1{\vphantom{#1}^{\#}}

\def\Z{\mathbb Z}
\def\R{\mathbb R}
\def\N{\mathbb N}
\def\cS{\mathcal S}
\def\e{{\bf e}}
\def\cP{\mathcal P}
\def\cQ{\mathcal Q}
\def\cT{\mathcal T}
\def\ord{{\rm ord\/}}
\def\eps{\varepsilon}
\def\ep{{\mathbf{\,e}}_p}
\def\epp{{\mathbf{\,e}}_{p-1}}

\def\xxx{\vskip5pt\hrule\vskip5pt}
\def\yyy{\vskip5pt\hrule\vskip2pt\hrule\vskip5pt}

\newcommand{\comm}[1]{\marginpar {\fbox{#1}}}
\newcommand{\1}{1\!{\rm l}}
\newcommand{\croix}{\overset\frown}


\title{\bf Almost All Palindromes\break
Are Composite}

\author{
        {\sc{William D.~Banks}\footnote{Corresponding author}} \\
        {\normalsize Dept.\ of Mathematics, University of Missouri} \\
        {\normalsize Columbia, MO 65211, USA} \\
        {\normalsize \tt bbanks@math.missouri.edu} \\
       \and
        {\sc{Derrick N.~Hart}} \\
        {\normalsize Dept.\ of Mathematics, University of Missouri} \\
        {\normalsize Columbia, MO 65211, USA} \\
        {\normalsize \tt hart@math.missouri.edu} \\
       \and
        {\sc{Mayumi Sakata}} \\
        {\normalsize Dept.\ of Mathematics, University of Missouri} \\
        {\normalsize Columbia, MO 65211, USA} \\
        {\normalsize \tt sakata@math.missouri.edu}}

\date{\today}

\maketitle

\vskip-15pt

\centerline{\it Never odd or even...}

\begin{abstract}
We study the distribution of palindromic numbers (with respect
to a fixed base $g\ge 2$) over certain congruence classes, and
we derive a nontrivial upper bound for the number of prime
palindromes $n\le x$ as $x\to\infty$.  Our results show that
almost all palindromes in a given base are composite.
\end{abstract}

\section{Introduction}

Fix once and for all an integer $g\ge 2$, and consider the
{\it base $g$ representation\/} of an arbitrary natural number $n\in\N$:
$$
n=\sum_{k=0}^{L-1}a_k(n)g^k.
$$
Here $a_k(n)\in\{0,1,\ldots,g-1\}$ for each $k=0,1,\ldots,L-1$,
and we assume that the leading digit $a_{L-1}(n)$ is nonzero.
The integer $n$ is said to be a {\it palindrome\/} if its digits satisfy
the symmetry condition:
$$
a_k(n)=a_{L-1-k}(n),\qquad k=0,1,\ldots,L-1.
$$
Let $\cP\subset\N$ denote the set of palindromes (in base $g$),
and for every positive real number $x$, let
$$
\cP(x)=\{n\le x\,|\,n\in\cP\}.
$$

In this paper, we study the distribution of palindromes in congruence
classes.  Using estimates for twisted Kloosterman sums to bound
exponential sums over the set $\cP_L$ of palindromes with
precisely $L$ digits, we show that the set $\cP(x)$ becomes
uniformly distributed (as $x\to\infty$) over the congruence
classes modulo $p$, where $p>g$ is any prime number for which
the multiplicative order ${\rm ord\/}_p(g)$ of $g$ in the
group $(\Z/p\Z)^\times$ is at least $3p^{1/2}$; see Corollary~\ref{cor:p}
to Theorem~\ref{thm:L to x} for a precise statement of this result.
We remark that, thanks to the work of Pappalardi~\cite{Pa}, almost all
primes $p$ satisfy the stronger condition
${\rm ord\/}_p(g)\ge p^{1/2}\exp((\log p)^c)$ where
$c$ is any constant less than $(1-\log 2)/2$;
see also~\cite{EM,IT}.

Using a variation of these techniques, we also show that the
set $\cP(x)$ becomes uniformly distributed (as $x\to\infty$)
over the congruence classes modulo $q$, where $q\ge 2$ is any integer
relatively prime to $g(g^2-1)$; see Corollary~\ref{cor:q}.
This latter result, although weaker than that obtained for primes
$p$ satisfying the condition $\ord_p(g)\ge 3p^{1/2}$,
allows us to deduce the main
result of this paper: {\it almost all palindromes in a given base
are composite\/}. More precisely, in Theorem~\ref{thm:main}, we show
that
$$
\card{\big\{}
\big\{n\in\cP(x)\,|\,\hbox{\rm $n$ is prime}\big\}
=O\(\card{\cP}\cP(x)\,\frac{\log\log\log x}{\log\log x}\),
\qquad x\to\infty,
$$
where the implied constant depends only on the base $g$.
This result appears to be the first of its kind in the literature.

{\bf Acknowledgments.} The authors would like to thank Florian Luca and
Igor Shparlinski, whose valuable observations on the original
manuscript led to significant improvements in our estimates.
During the preparation of this paper, W.~B.\ was supported in part
by NSF grant DMS-0070628.

\section{Preliminary Estimates}

For any integer $q\ge 2$, let $e_q(x)$ denote
the exponential function $\exp(2\pi ix/q)$, which is defined for
all $x\in\R$.  For any integer $c$ that is relatively prime to $q$,
let $\overline c$ denote an arbitrary multiplicative inverse
for $c$ modulo $q$; that is, $c\,\overline c\equiv 1\pmod q$.
Finally, let $d(q)$ be the number of positive integral divisors of $q$,
and let ${\rm ord\/}_q(g)$ be the smallest integer $t\ge 1$ such that
$g^t\equiv 1\pmod q$.

\begin{lem}
\label{lem:Kloost}
For all $q\ge 2$ with $\gcd(q,g)=1$ and all $a,b\in\Z$, we have
$$
\left|\sum_{k=1}^{\ord_q(g)}
e_q(ag^k+b\overline g^{\,k})\right|\le d(q)q^{1/2}\gcd(a,b,q)^{1/2}.
$$
\end{lem}

\begin{proof}
Consider the twisted Kloosterman sum
$$
K_\chi(a,b;q)=\sum_{\substack{1\le c\le q\\\gcd(c,q)=1}}
\chi(c)\,e_q(ac+b\overline c),
$$
where $\chi$ is a Dirichlet character modulo $q$.
The Weil-Estermann bound
$$
\big|K_\chi(a,b;q)\big|\le d(q)q^{1/2}\gcd(a,b,q)^{1/2}
$$
holds for such sums (the original proofs
by Weil~\cite{We} and Estermann~\cite{Es}
carry over for twisted sums with only slight modifications).
Averaging over all Dirichlet characters $\chi$ modulo $q$ for which
$\chi(g)=1$, it follows that
$$
\frac{\ord_q(g)}{\varphi(q)}\sum_\chi K_\chi(a,b;q)=
\sum_{\substack{1\le c\le q\\c\equiv g^k\pmod q,~\exists k}}
e_q(ac+b\overline c)
=\sum_{k=1}^{\ord_q(g)} e_q(ag^k+b\overline g^{\,k}).
$$
The result follows.
\end{proof}

\begin{lem}
\label{lem:sum roots}
The following bound holds for all $q\ge 2$, $k\ge 2$ and $h\in\Z$
provided that $q\!\not\vert~h$:
$$
\left|\sum_{a=0}^{k-1}e_q(ha)\right|
\le k\,\exp\(-\frac{4\gcd(h,q)^2}{q^2}\).
$$
\end{lem}

\begin{proof}
Let us write
$$
s(q,k,h)=\left|\sum_{a=0}^{k-1}e_q(ha)\right|.
$$
If $d=\gcd(h,q)$, then $s(q,k,h)=s(q/d,k,h/d)$,
hence it suffices to prove the assertion for the special case
where $\gcd(h,q)=1$, which we now assume.

Without loss of generality, we may also suppose
that $k\le q$.  Indeed, if $k\ge q+1$, then
we can express $k=mq+r$ with $0\le r\le q-1$ and simply observe that
$$
s(q,k,h)=s(q,r,h)\le r\le q-1
\le (q+1)\exp(-4/q^2) \le k\,\exp(-4/q^2).
$$

If $\gcd(h,q)=1$ and $2\le k\le q$, we have
\begin{eqnarray*}
s(q,k,h)^2=\sum_{a,b=0}^{k-1}e_q\(h(a-b)\)
&=&k+\sum_{\substack{a,b=0\\a\ne b}}^{k-1}
\cos\(\frac{2\pi h(a-b)}{q}\)\\
&\le&k+k(k-1)\cos(2\pi/q);
\end{eqnarray*}
therefore,
$$
\frac{s(q,k,h)^2}{k^2}\le 
\frac1k+\(1-\frac1k\)\cos(2\pi/q)
\le\frac12\(1+\cos(2\pi/q)\).
$$
Using the fact that
$1+\cos x\le 2\exp(-x^2/4)$ for $0\le x\le\pi$,
we obtain the desired result.
\end{proof}

\section{Exponential Sums over Palindromes}

For every $L\ge 1$, let $\cP_L$ denote the set of palindromes
(in base $g$) with precisely $L$ digits; that is,
$$
\cP_L=\{n\in\cP\,|\,g^{L-1}\le n<g^L\}.
$$

\begin{lem}
\label{lem:exp sum I}
Let $q\ge 2$ be an integer such that $p>g$ for every prime divisor
$p$ of $q$.  Then for every $c\in\Z$ such that
$$
\ord_q(g)>d(q)q^{1/2}\gcd(c,q)^{1/2},
$$
the exponential sum
$$
S_L(c)=\sum_{n\in\cP_L}e_q(cn)
$$
satisfies the bound
$$
\big|S_L(c)\big|\le\card{\cP}\cP_L\cdot\Theta_c^{(L-2\ord_q(g)-1)/4},
$$
where
$$
\Theta_c=\frac1g+\frac{(g-1)\,d(q)q^{1/2}\gcd(c,q)^{1/2}}{g\,\ord_q(g)}.
$$
\end{lem}

\begin{proof}
Since
\begin{eqnarray*}
S_{2L}(c)
&=&\sum_{a_0=1}^{g-1}\sum_{a_1=0}^{g-1}\ldots\sum_{a_{L-1}=0}^{g-1}
e_q\(\sum_{k=0}^{L-1}ca_k\(g^k+g^{2L-1-k}\)\)\\
&=&
\sum_{a_0=1}^{g-1}e_q\(ca_0\(1+g^{2L-1}\)\)
\prod_{k=1}^{L-1}\sum_{a_k=0}^{g-1}e_q\(ca_k\(g^k+g^{2L-1-k}\)\)
\end{eqnarray*}
and
\begin{eqnarray*}
\lefteqn{S_{2L+1}(c)
=\sum_{a_0=1}^{g-1}\sum_{a_1=0}^{g-1}\ldots\sum_{a_L=0}^{g-1}
e_q\(ca_Lg^L+\sum_{k=0}^{L-1}ca_k\(g^k+g^{2L-k}\)\)}\\
&&=
\sum_{a_0=1}^{g-1}e_q\(ca_0(1+g^{2L})\)
\sum_{a_L=0}^{g-1}e_q(ca_Lg^L)
\prod_{k=1}^{L-1}\sum_{a_k=0}^{g-1}e_q\(ca_k\(g^k+g^{2L-k}\)\),
\end{eqnarray*}
it follows that
$$
\big|S_{2L+\delta}(c)\big|\le (g-1)g^\delta
\prod_{k=1}^{L-1}\left|\sum_{a=0}^{g-1}
e_q\(ca\(g^k+g^{2L+\delta-1-k}\)\)\right|
$$
for all $L\ge 1$ and $\delta=0$ or $1$.

Put $N=\ord_q(g)$, and write $L-1=Nm+\ell$, where $m=\fl{(L-1)/N}$ and
$0\le\ell<N$.  Then, using the arithmetic-geometric mean
inequality, we derive that
\begin{eqnarray*}
\big|S_{2L+\delta}(c)\big|^2
&\le&(g-1)^2g^{2\ell+2\delta}
\prod_{k=1}^{Nm}\left|\sum_{a=0}^{g-1}
e_q\(ca\(g^k+g^{2L+\delta-1-k}\)\)\right|^2\\
&\le&(g-1)^2g^{2\ell+2\delta}\(\frac1{Nm}\sum_{k=1}^{Nm}
\left|\sum_{a=0}^{g-1}e_q\(ca\(g^k+g^{2L+\delta-1-k}\)\)
\right|^2\)^{\!\!\!Nm}\\
&=&(g-1)^2g^{2\ell+2\delta}\(\frac{T}{Nm}\)^{\!\!\!Nm},
\end{eqnarray*}
where
\begin{eqnarray*}
T&=&\sum_{k=1}^{Nm} \sum_{a,b=0}^{g-1}
e_q\(c(a-b)\(g^k+g^{2L+\delta-1-k}\)\)\\
&=&gNm+m\sum_{\substack{a,b=0\\a\ne b}}^{g-1}\sum_{k=1}^N
e_q\(c(a-b)\(g^k+g^{2L+\delta-1-k}\)\).
\end{eqnarray*}
Since $p>g$ for every prime $p\,|\,q$, it follows that
$\gcd(a-b,q)=1$ whenever $a\ne b$.  Using Lemma~\ref{lem:Kloost},
we therefore obtain that
$$
\big|T\big|\le gNm+(g-1)gm\,d(q)q^{1/2}\gcd(c,q)^{1/2}.
$$
Consequently,
\begin{eqnarray*}
\big|S_{2L+\delta}(c)\big|^2
&\le&(g-1)^2g^{2\ell+2\delta}
\(\frac{gNm+(g-1)gm\,d(q)q^{1/2}\gcd(c,q)^{1/2}}{Nm}\)^{\!\!\!Nm}\\
&=&(g-1)^2g^{2Nm+2\ell+2\delta}
\(\frac{N+(g-1)\,d(q)q^{1/2}\gcd(c,q)^{1/2}}{gN}\)^{\!\!\!Nm}\\
&=&(g-1)^2g^{2L+2\delta-2}\Theta_c^{L-\ell-1}.
\end{eqnarray*}
Since $\card{\cP}\cP_{2L+\delta}=(g-1)g^{L+\delta-1}$,
the result follows.
\end{proof}

\newpage

\begin{lem}
\label{lem:exp sum II}
Let $q\ge 2$ be an integer such that $\gcd\big(q,g(g^2-1)\big)=1$.
Then for every $c\in\Z$ such that $q\!\not\vert~c$, the exponential sum
$$
S_L(c)=\sum_{n\in\cP_L}e_q(cn)
$$
satisfies the bound
$$
\big|S_L(c)\big|\le\card{\cP}\cP_L\cdot
\exp\(-\frac{(L-5)\gcd(c,q)^2}{q^2}\).
$$
\end{lem}

\begin{proof}
As in the proof of Lemma~\ref{lem:exp sum I}, we have
$$
\big|S_{2L+\delta}(c)\big|\le (g-1)g^\delta
\prod_{k=1}^{L-1}\left|\sum_{a=0}^{g-1}
e_q\(ca\(g^k+g^{2L+\delta-1-k}\)\)\right|
$$
for all $L\ge 1$ and $\delta=0$ or $1$. Let
\begin{eqnarray*}
B&=&\{1\le k\le L-1\,|\,\hbox{$q$ divides
$c(g^k+g^{2L+\delta-1-k})$}\},\\
G&=&\{1\le k\le L-1\,|\,\hbox{$q$ does not divide
$c(g^k+g^{2L+\delta-1-k})$}\}.
\end{eqnarray*}
Using Lemma~\ref{lem:sum roots} to estimate individual terms in
the preceding product when $k\in G$, and using the trivial estimate
when $k\in B$, we obtain that
\begin{eqnarray*}
\big|S_{2L+\delta}(c)\big|
&\le&(g-1)g^{(\delta+\card{G}G+\card{B}\!B)}
\exp\(-\frac{4\gcd(c,q)^2\,\card{G}G}{q^2}\)\\
&=&\card{\cP}\cP_{2L+\delta}\cdot
\exp\(-\frac{4\gcd(c,q)^2\,\card{G}G}{q^2}\).
\end{eqnarray*}
Now let $f=q/\gcd(c,q)$.  Since $q$ does not divide $c$,
we have $f\ge 2$, and the stated condition on $q$ implies that
$\ord_f(g^2)\ge 2$.  Thus, if $k$ and $\ell$ both lie in $B$, then
$$
(g^2)^k\equiv g^{2L+\delta-1}\equiv (g^2)^\ell\pmod f,\qquad
k\equiv\ell\pmod{\ord_f(g^2)}.
$$
We therefore see that
\begin{eqnarray*}
\card{B}\!B&\le&1+\fl{(L-2)/2}=\fl{L/2},\\
\card{G}G&\ge&L-1-\fl{L/2}\ge L/2-1\ge (2L+\delta-5)/4,
\end{eqnarray*}
and the result follows.
\end{proof}

\section{Distribution of Palindromes}

\begin{prop}
\label{prop:P_L p}
Let $p>g$ be a prime number such that $\ord_p(g)\ge 3p^{1/2}$.
Then for every $L\ge 10p-5$, the following estimate holds
for all $a\in\Z$:
$$
\left|\card{\big\{}
\big\{n\in\cP_L\,|\,n\equiv a\pmod p\big\}-\frac{\card{\cP}\cP_L}p\right|
<\frac{\card{\cP}\cP_L}p\,(0.99)^L.
$$
\end{prop}

\begin{proof}
Using the relation
$$
\frac1p\sum_{c=0}^{p-1}e_p(cm)
=\left\{\begin{array}{ll}
1&\qquad\text{if $m\equiv 0\pmod p$},\\
0&\qquad\text{otherwise},\end{array} \right.
$$
it follows that
\begin{eqnarray*}
\card{\big\{}\big\{n\in\cP_L\,|\,n\equiv a\pmod p\big\}
&=&\sum_{n\in\cP_L}\frac1p\sum_{c=0}^{p-1}e_p\(c(n-a)\)\\
&=&\frac1p\sum_{c=0}^{p-1}e_p(-ca)\sum_{n\in\cP_L} e_p(cn)\\
&=&\frac{\card{\cP}\cP_L}p+\frac1p\sum_{c=1}^{p-1}e_p(-ca)S_L(c),
\end{eqnarray*}
where $S_L(c)$ is the exponential sum
considered in Lemma~\ref{lem:exp sum I}.  Therefore
\begin{eqnarray*}
\left|\card{\big\{}
\big\{n\in\cP_L\,|\,n\equiv a\pmod p\big\}-\frac{\card{\cP}\cP_L}p\right|
&\le& \frac{1}{p} \sum_{c=1}^{p-1} \big|S_L(c)\big| \\
&\le&\frac{\card{\cP}\cP_L}p\sum_{c=1}^{p-1}\Theta_c^{(L-2\,\ord_p(g)-1)/4},
\end{eqnarray*}
where for each $c=1,\ldots,p-1$, we have
$$
\Theta_c=\frac1g+\frac{2(g-1)p^{1/2}}{g\,\ord_p(g)}
\le\frac1g+\frac{2(g-1)}{3g}=\frac{2g+1}{3g}\le\frac56
$$
since $g\ge 2$.  Also,
$$
(L-2\,\ord_p(g)-1)/4\ge (L-2p+1)/4\ge L/5,\qquad L\ge 10p-5.
$$
Consequently,
$$
\left|\card{\big\{}
\big\{n\in\cP_L\,|\,n\equiv a\pmod p\big\}-\frac{\card{\cP}\cP_L}p\right|
\le\frac{\card{\cP}\cP_L}p\,(p-1)\(\frac56\)^{L/5}.
$$
Finally, remarking that the condition $\ord_p(g)\ge 3p^{1/2}$ implies
that $p\ge 11$, we have
$$
(p-1)\(\frac56\)^{L/5}<(0.99)^L,\qquad L\ge 10p-5.
$$
This completes the proof.
\end{proof}

\begin{prop}
\label{prop:P_L q}
Let $q\ge 2$ be an integer such that $\gcd\big(q,g(g^2-1)\big)=1$.
Then for every $L\ge 10+2q^2\log q$,
the following estimate holds for all $a\in\Z$:
$$
\left|\card{\big\{}
\big\{n\in\cP_L\,|\,n\equiv a\pmod q\big\}-\frac{\card{\cP}\cP_L}q\right|
<\frac{\card{\cP}\cP_L}{q}\,\exp\(-\frac{L}{2q^2}\).
$$
\end{prop}

\begin{proof}
Using the relation
$$
\frac1q\sum_{c=0}^{q-1}e_q(cm)
=\left\{\begin{array}{ll}
1&\qquad\text{if $m\equiv 0\pmod q$},\\
0&\qquad\text{otherwise},\end{array} \right.
$$
it follows that
\begin{eqnarray*}
\card{\big\{}\big\{n\in\cP_L\,|\,n\equiv a\pmod q\big\}
&=&\sum_{n\in\cP_L}\frac1q\sum_{c=0}^{q-1}e_q\(c(n-a)\)\\
&=&\frac1q\sum_{c=0}^{q-1}e_q(-ca)\sum_{n\in\cP_L} e_q(cn)\\
&=&\frac{\card{\cP}\cP_L}q+\frac1q\sum_{c=1}^{q-1}e_q(-ca)S_L(c),
\end{eqnarray*}
where $S_L(c)$ is the exponential sum
considered in Lemma~\ref{lem:exp sum II}.
If $1\le c\le q-1$, then $q\!\not\vert~c$, hence
by Lemma~\ref{lem:exp sum II} we derive the estimate:
\begin{eqnarray*}
\big|S_L(c)\big|
&\le&\frac{\card{\cP}\cP_L}{q}\,
\exp\(\log q-\frac{(L-5)\gcd(c,q)^2}{q^2}\)\\
&\le&\frac{\card{\cP}\cP_L}{q}\,\exp\(\log q-\frac{L-5}{q^2}\)
\le\frac{\card{\cP}\cP_L}{q}\,\exp\(-\frac{L}{2q^2}\),
\end{eqnarray*}
the last inequality following from the stated condition on $L$.
The result follows immediately.
\end{proof}

\begin{thm}
\label{thm:L to x}
Let $q\ge 2$ be a fixed integer, and
suppose that there exist constants $A\ge 1$ and $\sqrt{2/3}\le\xi<1$,
depending only on $q$, such that
$$
\left|\card{\big\{}
\big\{n\in\cP_L\,|\,n\equiv a\pmod q\big\}-\frac{\card{\cP}\cP_L}q\right|
\le \card{\cP}\cP_L\cdot A\,\xi^L
$$
for all $L\ge 1$ and $a\in\Z$.  Then for some constant $B\ge 1$
that depends only on $g$, the following estimate holds for
all $x\ge 1$ and $a\in\Z$:
$$
\left|\card{\big\{}
\big\{n\in\cP(x)\,|\,n\equiv a\pmod q\big\}-
\frac{\card{\cP}\cP(x)}q\right|
\le \card{\cP}\cP(x)\cdot AB\,\xi^{(\log x)/(2\log g)}.
$$
\end{thm}

\begin{proof} We remark that the condition $\xi\ge\sqrt{2/3}$
guarantees that $g\xi^2$ is bounded below by an absolute constant greater
than $1$; since $g\ge 2$, we have
$$
\frac{g-1}{g\xi^2-1}\le\frac{g-1}{\tfrac23 g-1}\le 3.
$$
For all $L\ge 1$, $x\ge y>0$, and $a\in\Z$, let us denote
\begin{eqnarray*}
\cP_a&=&\{n\in\cP\,|\,n\equiv a\pmod q\},\\
\cP_{a,L}&=&\{n\in\cP_a\,|\,g^{L-1}\le n<g^L\},\\
\cP_a(x)&=&\{n\in\cP_a\,|\,n\le x\},\\
\cP_a(y;x)&=&\{n\in\cP_a\,|\,y<n\le x\}.
\end{eqnarray*}
We also denote
$$
\cP(y;x)=\{n\in\cP\,|\,y<n\le x\}.
$$
In what follows, the implied constants in the symbol ``$O$''
may depend on $g$ but are absolute otherwise.
We recall that the notation $U=O(V)$ for positive functions $U$ and $V$
is equivalent to $U\le cV$ for some constant $c$.

Let $a\in\Z$ be fixed in what follows,
and suppose that $g^{2M+\delta-1}\le x<g^{2M+\delta}$, where $M$
is an integer and $\delta=0$ or $1$.  We observe that
\begin{equation}
\label{eq:P(x) equals}
\card{\cP}\cP(x)=\card{\cP}\cP(g^{2M+\delta-1})
+\card{\cP}\cP(g^{2M+\delta-1};x),
\end{equation}
and that
$$
\card{\cP}\cP_a(x)=\card{\cP}\cP_a(g^{2M+\delta-1})
+\card{\cP}\cP_a(g^{2M+\delta-1};x).
$$
Our goal is to estimate
\begin{eqnarray}
\label{eq:P(x) est}
\lefteqn{\left|\card{\cP}\cP_a(x)-\frac{\card{\cP}\cP(x)}q\right|}\\
&&\hskip-15pt\le\left|\card{\cP}\cP_a(g^{2M+\delta-1})
-\frac{\card{\cP}\cP(g^{2M+\delta-1})}q\right|
+\left|\card{\cP}\cP_a(g^{2M+\delta-1};x)
-\frac{\card{\cP}\cP(g^{2M+\delta-1};x)}q\right|.\nonumber
\end{eqnarray}

Since the integer $g^{2M+\delta-1}$
is {\it not\/} a palindrome (a fact that is only used to simplify our
notation), we have by a straightforward calculation:
\begin{equation}
\label{eq:P(M) calc}
\card{\cP}\cP(g^{2M+\delta-1})=g^M+g^{M+\delta-1}-2.
\end{equation}
On the other hand,
\begin{eqnarray*}
\lefteqn{\card{\cP}\cP_a(g^{2M+\delta-1})
=\sum_{L=1}^{2M+\delta-1}\card{\cP}\cP_{a,L}
=\sum_{\ell=0}^{M-1}\card{\cP}\cP_{a,2\ell+1}
+\sum_{\ell=1}^{M+\delta-1}\card{\cP}\cP_{a,2\ell}}\\
&&\hskip-8pt=\sum_{\ell=0}^{M-1}\(\card{\cP}\cP_{a,2\ell+1}-
\frac{\card{\cP}\cP_{2\ell+1}}q+\frac{\card{\cP}\cP_{2\ell+1}}q\)
+\sum_{\ell=1}^{M+\delta-1}\(\card{\cP}\cP_{a,2\ell}-
\frac{\card{\cP}\cP_{2\ell}}q+\frac{\card{\cP}\cP_{2\ell}}q\)\\
&&\hskip-8pt=\frac{\card{\cP}\cP(g^{2M+\delta-1})}q+
\sum_{\ell=0}^{M-1}\(\card{\cP}\cP_{a,2\ell+1}-
\frac{\card{\cP}\cP_{2\ell+1}}q\)
+\sum_{\ell=1}^{M+\delta-1}\(\card{\cP}\cP_{a,2\ell}-
\frac{\card{\cP}\cP_{2\ell}}q\).
\end{eqnarray*}
Using the hypothesis of the theorem, it therefore follows that
$$
\left|\card{\cP}\cP_a(g^{2M+\delta-1})
-\frac{\card{\cP}\cP(g^{2M+\delta-1})}q\right|
\le\sum_{\ell=0}^{M-1}\card{\cP}\cP_{2\ell+1}\cdot A\,\xi^{2\ell+1}
+\sum_{\ell=1}^{M+\delta-1}\card{\cP}\cP_{2\ell}\cdot A\,\xi^{2\ell}.
$$
Since
\begin{eqnarray*}
\lefteqn{\sum_{\ell=0}^{M-1}
\card{\cP}\cP_{2\ell+1}\,\xi^{2\ell+1}
+\sum_{\ell=1}^{M+\delta-1}
\card{\cP}\cP_{2\ell}\,\xi^{2\ell}}\\
&&\qquad\qquad=\sum_{\ell=0}^{M-1}(g-1)g^\ell\xi^{2\ell+1}
+\sum_{\ell=1}^{M+\delta-1}(g-1)g^{\ell-1}\xi^{2\ell}\\
&&\qquad\qquad<\frac{g-1}{g\xi^2-1}\,
\(g^M\xi^{2M+1}+g^{M+\delta-1}\xi^{2M+2\delta}\)
=O\(g^M\xi^{2M}\),
\end{eqnarray*}
we see that
\begin{equation}
\label{eq:Pa(M)}
\left|\card{\cP}\cP_a(g^{2M+\delta-1})
-\frac{\card{\cP}\cP(g^{2M+\delta-1})}q\right|=O\(Ag^M\xi^{2M}\).
\end{equation}

We now turn to the more delicate estimation of 
$\card{\cP}\cP_a(g^{2M+\delta-1};x)$.  To this end, put
$M=K+L$, where $K$ and $L$ are positive integers to be selected
later.  Examining the base $g$ representation of an arbitrary
palindrome $n$ in $\cP_{2M+\delta}$,
we see that $n$ may be expressed either in the form
$$
n=n_1+g^{K+\mu}n_2+g^{K+2L+\delta}n_3,
$$
or the form
$$
n=n_1+g^{K+2L+\delta}n_3,
$$
where
\begin{equation}
\label{eq:comp}
1\le n_1<g^K,\qquad
g^{K-1}\le n_3<g^K,\qquad n_1+g^Kn_3\in\cP_{2K},
\end{equation}
and, in the former case, $n_2\in\cP_{2L+\delta-2\mu}$ for some
$0\le\mu\le L+\delta-1$.  The integers $n_1,n_2,n_3,\mu$
are uniquely determined by $n$.  We call $n_3$ the
{\it K-signature\/} of $n$ and write $s_K(n)=n_3$.
The integer $n_1$ is uniquely determined by $n_3$
together with the first and third conditions of~(\ref{eq:comp});
we call $n_1$ the {\it K-complement\/} of $n_3$ and write
$c_K(n_3)=n_1$.

Note that the number of palindromes $n\in\cP_{2M+\delta}$
with a fixed $K$-signature $s_K(n)=n_3$ is precisely
\begin{equation}
\label{eq:pal sign}
1+\sum_{\mu=0}^{L+\delta-1}\card{\cP}\cP_{2L+\delta-2\mu}
=1+\sum_{\mu=0}^{L+\delta-1}(g-1)g^{L+\delta-\mu-1}=g^{L+\delta}.
\end{equation}

Now, given $x$ in the range $g^{2M+\delta-1}\le x<g^{2M+\delta}$,
let $y$ be the palindrome in $\cP_{2M+\delta}$ defined by
$$
y=y_1+g^K(g^{2L+\delta}-1)+g^{K+2L+\delta}y_3,
$$
where
$$y_3=\left\{\begin{array}{ll}
\fl{x/g^{K+2L+\delta}}+1&
\qquad\mbox{if $g^{2M+\delta-1}\le x<g^{2M+\delta-1/2}$},\\ \\
\fl{x/g^{K+2L+\delta}}-1&
\qquad\mbox{if $g^{2M+\delta-1/2}\le x<g^{2M+\delta}$},\end{array}\right.
$$
and $y_1=c_K(y_3)$.  If $x$ lies in the smaller range, then $x<y$,
while $y<x$ if $x$ lies in the larger range.  In either case, we have
\begin{equation}
\label{eq:P(M,x) calc}
\left|\card{\cP}\cP(g^{2M+\delta-1};x)-
\card{\cP}\cP(g^{2M+\delta-1};y)\right|=O(g^L)
\end{equation}
and
$$
\left|\card{\cP}\cP_a(g^{2M+\delta-1};x)-
\card{\cP}\cP_a(g^{2M+\delta-1};y)\right|
=O(g^L),
$$
since there are at most $O(1)$ distinct $K$-signatures for
palindromes between $x$ and $y$.  Consequently,
\begin{eqnarray}
\label{eq:Pa(M,x)-Pa(M,y)}
\lefteqn{\left|\card{\cP}\cP_a(g^{2M+\delta-1};x)
-\frac{\card{\cP}\cP(g^{2M+\delta-1};x)}q\right|}\\
&&\qquad\qquad\qquad\qquad=\left|\card{\cP}\cP_a(g^{2M+\delta-1};y)
-\frac{\card{\cP}\cP(g^{2M+\delta-1};y)}q\right|+O(g^L).\nonumber
\end{eqnarray}

Now, if $n\in\cP(g^{2M+\delta-1};y)$, then its $K$-signature lies in
the range
$$
g^{K-1}\le s_K(n)\le y_3.
$$
Thus,
\begin{equation}
\label{eq:P(m,y) calc}
\card{\cP}\cP(g^{2M+\delta-1};y)=(y_3-g^{K-1}+1)g^{L+\delta}.
\end{equation}
On the other hand, if $n\in\cP_a(g^{2M+\delta-1};y)$
with $s_K(n)=n_3$, then either
$$
n=n_1+g^{K+\mu}n_2+g^{K+2L+\delta}n_3\equiv a\pmod q
$$
or
$$
n=n_1+g^{K+2L+\delta}n_3\equiv a\pmod q,
$$
depending on the form of $n$.  In the latter case, there is
at most one such palindrome $n$ (for each fixed $K$-signature
$n_3$), while in the former case, since
$$
n_2\equiv g^{-K-\mu}\big(a-c_K(n_3)-g^{K+2L+\delta}n_3\big)\pmod q,
$$
the number of such palindromes $n$ is $\card{\cP}\cP_{b,2L+\delta-2\mu}$
for each $0\le\mu\le L+\delta-1$,
where
$$
b=b(n_3,\mu)=g^{-K-\mu}(a-c_K(n_3)-g^{K+2L+\delta}n_3).
$$
Hence, using~(\ref{eq:pal sign}), we derive that
\begin{eqnarray*}
\lefteqn{\card{\cP}\cP_a(g^{2M+\delta-1};y)
=\sum_{n_3=g^{K-1}}^{y_3}
\sum_{\mu=0}^{L+\delta-1}\card{\cP}\cP_{b,2L+\delta-2\mu}+O(g^K)}\\
&&\hskip-4pt=\sum_{n_3=g^{K-1}}^{y_3}
\(\frac1q+\sum_{\mu=0}^{L+\delta-1}
\frac{\card{\cP}\cP_{2L+\delta-2\mu}}q\)\\
&&\hskip-4pt\qquad\qquad\qquad+\sum_{n_3=g^{K-1}}^{y_3}
\sum_{\mu=0}^{L+\delta-1}\(\card{\cP}\cP_{b,2L+\delta-2\mu}
-\frac{\card{\cP}\cP_{2L+\delta-2\mu}}q\)+O(g^K)\\
&&\hskip-4pt=\frac{\card{\cP}\cP(g^{2M+\delta-1};y)}q
+\sum_{n_3=g^{K-1}}^{y_3}
\sum_{\mu=0}^{L+\delta-1}\(\card{\cP}\cP_{b,2L+\delta-2\mu}
-\frac{\card{\cP}\cP_{2L+\delta-2\mu}}q\)+O(g^K).
\end{eqnarray*}
Using the hypothesis of the theorem, it therefore follows that
\begin{eqnarray*}
\lefteqn{\left|\card{\cP}\cP_a(g^{2M+\delta-1};y)
-\frac{\card{\cP}\cP(g^{2M+\delta-1};y)}q\right|}\\
&&\qquad\qquad\le\sum_{n_3=g^{K-1}}^{y_3}
\sum_{\mu=0}^{L+\delta-1}\card{\cP}\cP_{2L+\delta-2\mu}
\cdot A\,\xi^{2L+\delta-2\mu}+O(g^K)\\
&&\qquad\qquad=\sum_{n_3=g^{K-1}}^{y_3}
\sum_{\mu=0}^{L+\delta-1}(g-1)g^{L+\delta-\mu-1}
\cdot A\,\xi^{2L+\delta-2\mu}+O(g^K)\\
&&\qquad\qquad<A(y_3-g^{K-1}+1)\(
\frac{g-1}{g\xi^2-1}\,g^{L+\delta}\xi^{2L+\delta+2}\)+O(g^K),
\end{eqnarray*}
and consequently,
$$
\left|\card{\cP}\cP_a(g^{2M+\delta-1};y)
-\frac{\card{\cP}\cP(g^{2M+\delta-1};y)}q\right|
=O\(Ag^M\xi^{2L}\)+O(g^K).
$$

Using this estimate together with~(\ref{eq:P(x) est}),
(\ref{eq:Pa(M)}) and~(\ref{eq:Pa(M,x)-Pa(M,y)}), it follows that
$$
\left|\card{\cP}\cP_a(x)-\frac{\card{\cP}\cP(x)}q\right|=
O\(Ag^M\xi^{2L}+g^L+g^K\).
$$
We now choose integers $K=M/2+O(1)$ and $L=M/2+O(1)$ such that $K+L=M$.
Since $g\xi^2>1$ and $A\ge 1$, we have
$$
\max\{g^K,g^L\}=O(g^{M/2})=O\(Ag^{M/2}(g\xi^2)^{M/2}\)
=O\(Ag^M\xi^M\),
$$
therefore
$$
\left|\card{\cP}\cP_a(x)-\frac{\card{\cP}\cP(x)}q\right|=
O\(Ag^M\xi^M\).
$$
To complete the proof, we need only observe that
$$
\xi^M=O\(\xi^{(\log x)/(2\log g)}\)
$$
for $x$ in the range $g^{2M+\delta-1}\le x<g^{2M+\delta}$,
and using (\ref{eq:P(x) equals}), (\ref{eq:P(M) calc}),
(\ref{eq:P(M,x) calc}) and (\ref{eq:P(m,y) calc})
together with our choice of $y_3$, it follows that
$$
\card{\cP}\cP(x)=g^M+\frac{x}{g^M}+O(g^{M/2});
$$
thus $g^M=O\(\card{\cP}\cP(x)\)$.
\end{proof}

Using Theorem~\ref{thm:L to x}, we can now derive two immediate corollaries.

\begin{cor}
\label{cor:p}
Let $p>g$ be a prime number such that $\ord_p(g)\ge 3p^{1/2}$.
Then for some constant $C>0$, depending only on $g$,
the following estimate holds for all $x\ge 1$ and $a\in\Z$:
$$
\left|\card{\big\{}
\big\{n\in\cP(x)\,|\,n\equiv a\pmod p\big\}-\frac{\card{\cP}\cP(x)}p\right|
\le\card{\cP}\cP(x)\cdot C\,(0.99)^{\frac{\log x}{2\log g}-10p}.
$$
\end{cor}

\begin{proof}
Using the trivial estimate
$$
\left|\card{\big\{}
\big\{n\in\cP_L\,|\,n\equiv a\pmod p\big\}-\frac{\card{\cP}\cP_L}p\right|
\le \card{\cP}\cP_L
$$
for $1\le L\le 10p-6$, it follows from Proposition~\ref{prop:P_L p} that
the estimate
$$
\left|\card{\big\{}
\big\{n\in\cP_L\,|\,n\equiv a\pmod p\big\}-\frac{\card{\cP}\cP_L}p\right|
\le\card{\cP}\cP_L\cdot(0.99)^{L-10p+6}
$$
holds for all $L\ge 1$ and $a\in\Z$.  The result now
follows immediately from Theorem~\ref{thm:L to x}.
\end{proof}

\begin{cor}
\label{cor:q}
Let $q\ge 2$ be an integer such that $\gcd\big(q,g(g^2-1)\big)=1$.
Then for some constant $C>0$, depending only on $g$,
the following estimate holds for all $x\ge 1$ and $a\in\Z$:
$$
\left|\card{\big\{}
\big\{n\in\cP(x)\,|\,n\equiv a\pmod q\big\}
-\frac{\card{\cP}\cP(x)}q\right|
\le\card{\cP}\cP(x)\cdot
C\,q\exp\(-\frac{\log x}{4q^2\log g}\).
$$
\end{cor}

\begin{proof}
Using the trivial estimate
$$
\left|\card{\big\{}
\big\{n\in\cP_L\,|\,n\equiv a\pmod q\big\}-\frac{\card{\cP}\cP_L}q\right|
\le \card{\cP}\cP_L
$$
for $1\le L<10+2q^2\log q$,
it follows from Proposition~\ref{prop:P_L q} that the estimate
$$
\left|\card{\big\{}
\big\{n\in\cP_L\,|\,n\equiv a\pmod q\big\}
-\frac{\card{\cP}\cP_L}q\right|
\le\card{\cP}\cP_L\,\exp\(-\frac{(L-10-2q^2\log q)}{2q^2}\)
$$
holds for all $L\ge 1$ and $a\in\Z$.  The result now
follows immediately from Theorem~\ref{thm:L to x}.
\end{proof}

\section{Prime Palindromes}

We now come to the main result of this paper.

\begin{thm}
\label{thm:main}
As $x\to\infty$, we have
$$
\card{\big\{}
\big\{n\in\cP(x)\,|\,\hbox{\rm $n$ is prime}\big\}
=O\(\card{\cP}\cP(x)\,\frac{\log\log\log x}{\log\log x}\),
$$
where the implied constant depends only on $g$.
\end{thm}

\begin{proof}
As in the proof of Theorem~\ref{thm:L to x}, all implied constants in
the symbol ``$O$'' may depend on $g$ but are absolute otherwise.

Assuming that $x$ is sufficiently large, let
$$
h=\fl{e\log\log\log x},\qquad
y=e^{-1}(\log x)^{1/4h}
=\exp\(\frac{\log\log x}{4e\log\log\log x}\)^{1+o(1)}.
$$
Let
$$
Q=Q(y)=\prod_{g^3<p\le y}p,
$$
where the product runs over prime numbers.
Note that $\gcd\big(Q,g(g^2-1)\big)=1$.
By Mertens' formula (see Theorem~11 in \S I.1.6 of~\cite{Te}),
we have the estimate
\begin{equation}
\label{eq:merten}
\frac{\varphi(Q)}Q=\prod_{g^3<p\le y}\(1-\frac1p\)
=O\((\log y)^{-1}\)
=O\(\frac{\log\log\log x}{\log\log x}\),
\end{equation}
where $\varphi(n)$ is the Euler function.

Now, if $n\in\cP(x)$ is prime,
either $\gcd(n,Q)=1$ or $n$ is a prime divisor of $Q$.
We apply Brun's combinatorial sieve in the form given by
Corollary~1.1 in~\S I.4.2 of~\cite{Te}:
$$
\card{\big\{}
\big\{n\in\cP(x)\,|\,\hbox{$n$ is prime}\big\}
\le y+\sum_{\substack{q\,|\,Q\\\omega(q)\le 2h}}
\mu(q)A_q,
$$
where $\mu(q)$ is the M\"obius function, $\omega(q)$ is
the number of distinct prime divisors of $q$, and
$$
A_q=\card{\{}\{n\in\cP(x)\,|\,n\equiv 0\pmod q\}.
$$
By Corollary~\ref{cor:q}, we see that
$$
A_q=\frac{\card{\cP}\cP(x)}q+
O\(\card{\cP}\cP(x)\,q\exp\(-\frac{\log x}{4q^2\log g}\)\).
$$
If $q\,|\,Q$ and $\omega(q)\le 2h$, then
$$
q\le y^{2h}=\frac{(\log x)^{1/2}}{e^{2h}},
$$
and since the number of such divisors $q$ is bounded by $y^{2h}$,
we have
\begin{eqnarray*}
\lefteqn{\sum_{\substack{q\,|\,Q\\\omega(q)\le 2h}}
q\exp\(-\frac{\log x}{4q^2\log g}\)
\le\frac{\log x}{e^{4h}}
\exp\(-\frac{e^{4h}}{4\log g}\)}\\
&&\qquad\qquad=\exp\(\log\log x-4h-\frac{e^{4h}}{4\log g}\)
=O\(\frac{1}{\log x}\),
\end{eqnarray*}
since $h=\fl{e\log\log\log x}$.  Therefore,
\begin{eqnarray*}
\lefteqn{\card{\big\{}
\big\{n\in\cP(x)\,|\,\hbox{$n$ is prime}\big\}}\\
&&\qquad\le y+\card{\cP}\cP(x)
\sum_{q\,|\,Q}\frac{\mu(q)}q+
O\(\card{\cP}\cP(x)
\sum_{\substack{q\,|\,Q\\\omega(q)>2h}}
\frac1q+O\(\frac{\card{\cP}\cP(x)}{\log x}\)\).
\end{eqnarray*}
Since $y=x^{o(1)}$ and $x^{1/2}=O\(\card{\cP}\cP(x)\)$, the first
term in this estimate is negligible.  Also, using~(\ref{eq:merten}),
we have
$$
\card{\cP}\cP(x)\sum_{q\,|\,Q}\frac{\mu(q)}q
=\card{\cP}\cP(x)\prod_{g^3<p\le y}\(1-\frac1p\)
=O\(\card{\cP}\cP(x)\,\frac{\log\log\log x}{\log\log x}\).
$$
Finally, we have
$$
\sum_{\substack{q\,|\,Q\\\omega(q)>2h}}\frac1q
\le\sum_{\substack{q\,|\,Q\\\omega(q)>2h}}\frac{e^{\omega(q)-2h}}q
\le e^{-2h}\prod_{p\le y}(1+e/p)
\le\exp\(-2h+e\sum_{p\le y}1/p\).
$$
Observing that
$$
\sum_{p\le y}\frac1p=(\log\log y)(1+o(1))
=(\log\log\log x)(1+o(1)),
$$
by our choice of $h$ it follows that
$$
\card{\cP}\cP(x)\hskip-3.5pt
\sum_{\substack{q\,|\,Q\\\omega(q)>2h}}\hskip-3.5pt
\frac1q\le
\card{\cP}\cP(x)\exp\((\log\log\log x)(-e+o(1))\)
=O\(\frac{\card{\cP}\cP(x)}{(\log\log x)^2}\).
$$
This completes the proof.
\end{proof}

\section{Remarks and Open Problems}

Using estimates from~\cite{CPR}, it is possible to establish a
version of Lemma~\ref{lem:exp sum I} in the case where
$q=p$ is prime with ${\rm ord\/}_p(g)\gg \log p$; this yields
analogues of Proposition~\ref{prop:P_L p} and Corollary~\ref{cor:p}
under the weaker assumption on ${\rm ord\/}_p(g)$, however the uniform
constant $0.99$ in those results must be replaced by a term like
$\exp(-(\log\log p)^{-c})$ for some constant $c>0$.

It seems natural to conjecture that the set of palindromes should
behave as ``random'' integers, thus one might expect that the
asymptotic relation
$$
\card{\big\{}
\big\{n\in\cP(x)\,|\,\hbox{$n$ is prime}\big\}
\sim C\,\frac{\card{\cP}\cP(x)}{\log x}
$$
holds for some constant $C>0$.  While this question seems out of
reach at the moment, it should be feasible to derive the upper bound
$$
\card{\big\{}
\big\{n\in\cP(x)\,|\,\hbox{$n$ is prime}\big\}
=O\(\frac{\card{\cP}\cP(x)}{\log x}\)
$$
using more sophisticated sieving techniques coupled with better estimates
for the distribution of palindromes in congruence classes.
It is still an open problem to show the existence
of infinitely many prime palindromes for any fixed base $g\ge 2$.

\end{document}